\documentclass[11pt,a4paper]{amsproc}
\usepackage[a4paper,left=22mm,right=22mm,top=22mm,bottom=24mm,
            headheight=14pt,headsep=6mm,footskip=10mm]{geometry}
\usepackage[T1]{fontenc}
\usepackage{amsmath,amssymb,amsthm,mathtools,microtype}
\usepackage[colorlinks=true,linkcolor=blue,citecolor=blue,urlcolor=blue]{hyperref}
\newtheorem{maintheorem}{Theorem}

\newtheorem{theorem}{Theorem}[section]
\newtheorem{lemma}[theorem]{Lemma}
\newtheorem{proposition}[theorem]{Proposition}

\theoremstyle{remark}

\newcommand{\R}{\mathbb R}
\newcommand{\C}{\mathbb C}
\newcommand{\M}{\mathcal M}
\newcommand{\one}{\mathbf 1}
\newcommand{\tr}{\tau}
\newcommand{\dd}{\,\mathrm d}
\newcommand{\supp}{\operatorname{supp}}
\newcommand{\spec}{\sigma}
\newcommand{\im}{\operatorname{Im}}
\newcommand{\re}{\operatorname{Re}}
\numberwithin{equation}{section}
\title[The BMV conjecture in finite von Neumann algebras]
{The Bessis--Moussa--Villani conjecture\break in finite von Neumann algebras}
\author{Liviu  P\u{a}unescu}
\address{Institute of Mathematics of the Romanian Academy, 21 Calea Grivitei, 010702 Bucharest Romania }
\email{liviu.paunescu@imar.ro}
\author{Florin R\u{a}dulescu}
\address{Dipartimento di Matematica, Universit\`a degli Studi di Roma
``Tor Vergata'', Via della Ricerca Scientifica 1, 00133 Roma, Italy}
\email{radulesc@mat.uniroma2.it}
\author{Andreas Thom}
\address{Institut f\"ur Geometrie, Technische Universit\"at Dresden,
Zellescher Weg 12--14, 01069 Dresden, Germany}
\email{andreas.thom@tu-dresden.de}
\date{16 November 2026}
\hypersetup{pdftitle={The Bessis--Moussa--Villani conjecture in finite von Neumann algebras},pdfauthor={Liviu P\u{a}unescu, Florin R\u{a}dulescu and Andreas Thom},
  pdfkeywords={BMV conjecture, finite von Neumann algebras, tracial joint spectral measures, moment problems, Herglotz representation, complete monotonicity}}
\begin{document}
\begin{abstract}
We prove the Bessis--Moussa--Villani conjecture, which is known for matrix algebras, for bounded operators in a
finite von Neumann algebra with a faithful normal tracial state.
Following the approach of Hein\"avaara, for a pair of positive operators, we construct a tracial joint spectral measure on the
punctured plane whose moments are the traces of symmetrized words.
We construct a scalar measure by holomorphic spectral selection and a
Lyapunov identity. An integral identity for its moments then gives
positivity of two coordinate-weighted distributions. From these we
obtain the planar measure and the Laplace representation of the trace
exponential.
\end{abstract}
\maketitle

\section{Introduction}

The Bessis--Moussa--Villani conjecture arose in 1975 from the study of
variational approximations in quantum statistical mechanics \cite{BMV}.
It asserts that, for Hermitian matrices $H,V$ with $V\geq0$, the function
$t\mapsto\operatorname{Tr}(e^{H-tV})$ is the Laplace transform of a
positive measure on $[0,\infty)$. Lieb and Seiringer showed that this is
equivalent to nonnegativity of the coefficients of
$\operatorname{Tr}((A+tB)^m)$ for positive matrices $A,B$ and positive
integers $m$ \cite{LS}. Stahl proved the conjecture by constructing
the representing measure and establishing its positivity through the
complex analysis of algebraic eigenvalue functions \cite{Stahl}.

The main result of this paper is a proof of the corresponding statements for bounded operators in finite
von Neumann algebras.

Let $\M$ be a finite von Neumann algebra with identity $\one$, equipped
with a faithful normal tracial state $\tr$. For operators
$A,B\in\M$ and integers $0\leq k\leq m$, let $S_{m,k}(A,B)$ denote the
sum of all words of length $m$ containing precisely $k$ copies of $B$
and $m-k$ copies of $A$. Thus
\[
 (A+tB)^m=\sum_{k=0}^m t^k S_{m,k}(A,B),
 \qquad S_{0,0}(A,B)=\one.
\]
We prove the following theorem.

\begin{maintheorem}\label{thm:bmv}
Let $(\M,\tr)$ be a finite von Neumann algebra with a faithful normal
tracial state. For all positive $A,B\in\M$ and all integers
$0\leq k\leq m$, we have
\[
 \tr(S_{m,k}(A,B))\geq0.
\]
If $H,V\in\M$, $H=H^*$ and $V\geq0$, then there is a unique finite
positive Borel measure $\rho_{H,V}$ supported in
$[\min\spec(V),\max\spec(V)]$ such that
\begin{equation}\label{eq:bmv}
 \tr(e^{H-tV})=\int e^{-ts}\dd\rho_{H,V}(s)
 \qquad(t\geq0).
\end{equation}
Its total mass is $\tr(e^H)$. In particular,
$t\mapsto\tr(e^{H-tV})$ is completely monotone on $[0,\infty)$.
\end{maintheorem}

Here complete monotonicity means that
$(-1)^k\frac{\mathrm d^k}{\mathrm dt^k}\tr(e^{H-tV})\geq0$ for every
integer $k\geq0$ and every $t\geq0$.

There were several algebraic approaches to the
matrix conjecture. Hillar established descent results for the two letter
counts; in particular, validity in one degree implies validity in every
smaller degree \cite[Theorem~1.7 and Corollary~1.8]{Hillar}.
Klep and Schweighofer combined these results with representations by sums
of Hermitian squares modulo commutators to prove the conjecture through
degree thirteen \cite{KS}. Burgdorf obtained such representations for
the fourth coefficient in every degree \cite{Burgdorf}.
Finally, Collins, Dykema and Torres-Ayala established obstructions to this
certificate method for an infinite family of word sums \cite{CDT}.

The finite-trace version was considered by Bo\.{z}ejko in terms of
positive definiteness of the trace exponential; he proved this property
for generalized Gaussian pairs \cite{Boz}. Klep and Schweighofer
discussed the passage from matrix trace positivity to arbitrary
$\mathrm{II}_1$ factors in connection with Connes' embedding conjecture,
already for $S_{6,3}(X^2,Y^2)$ \cite[Section~6.2]{KS}.
More recently, Pascoe and Woerdeman formulated the shuffle-word-embedding
conjecture and observed that it would imply the BMV inequalities in
every $\mathrm{II}_1$ factor \cite[Section~3.1]{PW}.
Theorem~\ref{thm:bmv} finally establishes these inequalities for arbitrary
positive pairs in a finite von Neumann algebra with a faithful normal
tracial state.

Our proof uses the tracial joint spectral measure introduced by
Hein\"avaara \cite{Hei}. For a pair of Hermitian matrices, this positive
planar measure represents traces of functions of linear combinations
through a fixed one-variable integral transform. Its existence gives
the matrix BMV theorem as a consequence and also yields isometric
embeddings of two-dimensional real subspaces of Schatten classes into
classical $L^p$ spaces.
Dinh, Le and Vuong used related matrix-pencil measures to obtain
positive planar representations of spectral pinching defects
\cite{DLV}. For $A>0$, $B\geq0$ and $p\geq1$, Rey-Bellet obtained a
positive Stieltjes representation of order $p$ for
$\operatorname{Tr}((A+tB)^{-p})$ from Hein\"avaara's measure and
recovered Stahl's measure by a limiting argument \cite{ReyBellet}.
We establish the following planar representation for positive
operators in $\M$.

\begin{maintheorem}\label{thm:planar}
Let $A,B\in\M$ be positive and put
$K=[0,\|A\|]\times[0,\|B\|]$. There is a unique positive Borel measure
$\mu_{A,B}$ on $K\setminus\{(0,0)\}$ with finite first moment
\begin{equation}\label{eq:first-moment}
 \int(a+b)\dd\mu_{A,B}(a,b)=\frac{\tr(A+B)}2
\end{equation}
and the following property. For every Borel function $f:\R\to\C$
satisfying
\[
 f(0)=0,\qquad
 \int_{-R}^{R}\frac{|f(u)|}{|u|}\dd u<\infty
 \quad\text{for every }R>0,
\]
define
\[
 H_f(s)=\int_0^1\frac{1-u}{u}f(us)\dd u.
\]
Then, for all $x,y\in\R$,
\begin{equation}\label{eq:planar-marginal}
 \int f(ax+by)\dd\mu_{A,B}(a,b)=\tr(H_f(xA+yB)),
\end{equation}
where the integral on the left is absolutely convergent.
\end{maintheorem}

The measure is finite on every compact subset of the punctured plane.
The function $H_f$ is bounded on each compact interval: if
$0<|s|\leq R$, then
\[
 |H_f(s)|\leq\int_{-R}^{R}\frac{|f(v)|}{|v|}\dd v,
\]
and $H_f(0)=0$. Thus the right side of \eqref{eq:planar-marginal} is
defined by the bounded Borel functional calculus.
Taking $f(s)=s^m$ gives
\begin{equation}\label{eq:intro-moments}
 \int(ax+by)^m\dd\mu_{A,B}(a,b)
 =\frac{\tr((xA+yB)^m)}{m(m+1)}\qquad(m\geq1).
\end{equation}
Comparing coefficients, we obtain the first assertion of
Theorem~\ref{thm:bmv}. We then pass to the exponential formulation
following the argument of Lieb and Seiringer \cite{LS}. The details of
this passage, including the construction and support of the Laplace
measure, are given in Section~\ref{sec:bmv}.

We outline the construction of the planar measure. The arguments
described here are carried out in Sections~\ref{sec:fourier}--\ref{sec:measure};
the main positivity step is Proposition~\ref{prop:weighted}.
Following \cite{Hei}, we set
\[
 G_{A,B}(\xi,\eta)
 =\int_0^1\frac{1-u}{u}
       \bigl[\tr(e^{iu(\xi A+\eta B)})-1\bigr]\dd u
\]
and take its Fourier transform in the sense of distributions, with the
normalization specified in Section~\ref{sec:fourier}.
Lemma~\ref{lem:fourier} shows that the resulting distribution $D$ is
supported in $K$ and has the moments prescribed by \eqref{eq:intro-moments}.
The proof uses the entire extension of the displayed function to obtain
the support bound, and differentiation at the origin gives the moments.
We then prove the weighted inequalities
\begin{equation}\label{eq:intro-positivity}
 D(a|p(a,b)|^2)\geq0,\qquad D(b|p(a,b)|^2)\geq0
 \qquad(p\in\C[a,b]).
\end{equation}
The scalar measure needed for these inequalities is constructed in
Proposition~\ref{prop:scalar}. For an arbitrary selfadjoint pair
$S,J\in\M$, it gives a finite positive measure $\nu_{S,J}$ supported in
the interval $[\min\spec(J),\max\spec(J)]$.
We write $X_+=(|X|+X)/2$ for the positive part of a selfadjoint
operator $X$. The measure is characterized by
\[
 \int\frac{\dd\nu_{S,J}(b)}{x-b}
 =\tr\left(\left[(x\one-J)^{-1/2}S(x\one-J)^{-1/2}\right]_+\right)
 \quad(x>\max\spec(J)).
\]
The corresponding matrix measure was obtained by Dinh, Le and Vuong
\cite[Lemma~3.2]{DLV} from the sum of the eigenvalues of
$S(z\one-J)^{-1}$ in the lower half-plane.
For invertible $S$, we use the Riesz idempotent $E_-(z)$ selecting this
spectral part and consider the holomorphic function
\[
 m(z)=\tr\bigl(S(z\one-J)^{-1}E_-(z)\bigr),\qquad \im z>0.
\]
The Lyapunov identity in Lemma~\ref{lem:lyapunov} gives the sign of its
imaginary part. The Herglotz representation theorem then gives a positive
measure. Its support and mass follow from the real boundary values and
the asymptotic behaviour along the positive imaginary axis.
Proposition~\ref{prop:scalar} also
passes to arbitrary selfadjoint pairs by norm approximation and proves
the weak continuity needed in the next step.

In Section~\ref{sec:positivity}, for invertible $A\geq0$, we apply this
construction to the family
\[
 S_a=a^{-1}\one-A^{-1},\qquad
 J_a=aA^{-1/2}BA^{-1/2},\qquad 0<a<\|A\|.
\]
We integrate the moments of these measures in $a$ and compare them with
the mixed moments of the distribution. The invariance of the trace of
positive spectral projections under congruence, proved in
Lemma~\ref{lem:inertia}, allows us to rewrite the generating function as
an integral of spectral tails. A count of marked words then gives the
moment identity \eqref{eq:slice-moments}. Proposition~\ref{prop:weighted}
uses this identity to obtain the first inequality in
\eqref{eq:intro-positivity}. When both operators are invertible,
interchanging $A$ and $B$ gives the second inequality. Norm approximation
then gives both inequalities for arbitrary positive operators.

Section~\ref{sec:measure} turns these inequalities into the planar
measure. Lemma~\ref{lem:positive-distribution} shows that a compactly
supported distribution which is nonnegative on polynomial squares is
a positive measure. Applying it to $aD$ and $bD$, we obtain two finite
positive measures. Their coordinate marginals show that they have no
atom at the origin. Dividing their sum by $a+b$, we obtain
$\mu_{A,B}$. The compensated identity \eqref{eq:compensated} recovers the
original distribution from this measure and gives all its moments of
positive degree. Finally, comparison of finite measures obtained by
weighting with the square of a linear form gives the full Borel identity
in Theorem~\ref{thm:planar}, and polynomial density proves uniqueness.

For a Borel set $E\subseteq\R$, we write
$\one_E(X)$ for its spectral projection. We use the restriction of
$\tr$ on every corner $P\M P$; its value on the corner identity is
$\tr(P)$. Spectra and inverses of corner operators are taken in that
corner.

\section{The Fourier distribution and its moments}\label{sec:fourier}

Fix positive $A,B\in\M$, and write
$\alpha=\|A\|$, $\beta=\|B\|$, and
$K=[0,\alpha]\times[0,\beta]$. Our Fourier convention is
\[
 \widehat\varphi(\xi,\eta)
 =\int_{\R^2}e^{-i(\xi a+\eta b)}\varphi(a,b)\dd a\dd b.
\]
For a Schwartz function $\varphi$, put
\begin{equation}\label{eq:distribution}
 D(\varphi)=\frac1{(2\pi)^2}
        \int_{\R^2}G_{A,B}(\xi,\eta)\widehat\varphi(\xi,\eta)
                   \dd\xi\dd\eta.
\end{equation}
Since $\|e^{iuX}-\one\|\leq u\|X\|$ for selfadjoint $X$ and $u\geq0$,
we have
\begin{equation}\label{eq:real-growth}
 |G_{A,B}(\xi,\eta)|\leq\tfrac12(\alpha|\xi|+\beta|\eta|).
\end{equation}
It follows that \eqref{eq:distribution} defines a tempered
distribution. The identity
$G_{A,B}(-\xi,-\eta)=\overline{G_{A,B}(\xi,\eta)}$ shows that $D$ is real.

\begin{lemma}\label{lem:fourier}
The distribution $D$ is supported in $K$. For all real $\xi,\eta$,
\begin{equation}\label{eq:transform}
 D(e^{i(\xi a+\eta b)})=G_{A,B}(\xi,\eta).
\end{equation}
Moreover, $D(1)=0$, and for all real $x,y$ and integers $n\geq1$,
\begin{equation}\label{eq:directional-moments}
 D((xa+yb)^n)=\frac{\tr((xA+yB)^n)}{n(n+1)}.
\end{equation}
\end{lemma}
\begin{proof}
The exponential series gives an entire extension of $G_{A,B}$ to
$\C^2$. For $z=(z_1,z_2)\in\C^2$, set
$L=i(z_1A+z_2B)$ and
$h=\alpha(-\im z_1)_++\beta(-\im z_2)_+$. Then
$\re L\leq h\one$. Differentiating $\|e^{uL}v\|^2$ for a vector $v$
in a Hilbert space representation of $\M$, we obtain
$\|e^{uL}\|\leq e^{uh}$. Hence
\[
 \|e^{uL}-\one\|
 \leq\int_0^u\|L\|e^{vh}\dd v
 \leq u(\alpha|z_1|+\beta|z_2|)e^h\quad(0\leq u\leq1),
\]
and therefore
\begin{equation}\label{eq:entire-bound}
 |G_{A,B}(z_1,z_2)|
 \leq\tfrac12(\alpha|z_1|+\beta|z_2|)
       e^{\alpha(-\im z_1)_++\beta(-\im z_2)_+}.
\end{equation}

To determine the support, let $\varphi\in C_c^\infty(\R^2)$ be supported
in $a>\alpha$.
Choose $\delta>0$ such that $a\geq\alpha+\delta$ on its support.
For fixed $R>0$, we move the first Fourier contour to $\xi-iR$.
Both factors in \eqref{eq:distribution} are entire in this variable.
On a fixed horizontal strip, $G_{A,B}$ grows at most linearly in the
real frequencies, whereas $\widehat\varphi$ decays faster than any
inverse power. We can therefore integrate over a rectangle and let its
vertical sides tend to infinity. These bounds also justify integration
in $\eta$, and we obtain
\[
 D(\varphi)=\frac1{(2\pi)^2}
   \int G_{A,B}(\xi-iR,\eta)\widehat\varphi(\xi-iR,\eta)
      \dd\xi\dd\eta.
\]
Integration by parts, applied to $e^{-Ra}\varphi(a,b)$, gives
\[
 |\widehat\varphi(\xi-iR,\eta)|
 \leq C_N(1+R)^{2N}e^{-(\alpha+\delta)R}
                 (1+\xi^2+\eta^2)^{-N}.
\]
Combining this estimate with \eqref{eq:entire-bound}, we obtain, for
$N\geq3$, a bound for $|D(\varphi)|$ by a polynomial in $R$ times
$e^{-\delta R}$.
Letting $R\to\infty$ gives $D(\varphi)=0$.
For a test supported in $a<0$, we shift to $\xi+iR$ and use the decay
of $e^{Ra}\varphi$. The shifts in the second variable give the same
conclusion on $b>\beta$ and on $b<0$. These four half-planes cover
$\R^2\setminus K$, so a partition of unity proves $\supp D\subseteq K$.

A compactly supported distribution acts on every smooth function by
inserting a fixed cutoff equal to one near its support. This extension
agrees with its original action on Schwartz functions. Indeed, a
Schwartz function vanishing near the support can be approximated in
the Schwartz topology by multiplying it by expanding compact cutoffs.
The resulting compactly supported tests also vanish near the support
and are therefore annihilated by the distribution.
For $\theta=(\xi,\eta)\in\R^2$, let
$\varphi_R(a,b)=e^{i(\xi a+\eta b)}e^{-(a^2+b^2)/(2R^2)}$.
Its Fourier transform is
$2\pi R^2e^{-R^2|\zeta-\theta|^2/2}$. Consequently,
\[
 D(\varphi_R)
 =\frac1{2\pi}\int_{\R^2}
       G_{A,B}(\theta+w/R)e^{-|w|^2/2}\dd w
 \longrightarrow G_{A,B}(\theta),
\]
where dominated convergence follows from \eqref{eq:real-growth}.
On a neighborhood of $K$, the functions $\varphi_R$ converge with
all derivatives to $e^{i(\xi a+\eta b)}$. This proves
\eqref{eq:transform}, and evaluation at $\theta=0$ gives $D(1)=0$.

Finally, put $X=xA+yB$. To obtain the moments, we differentiate
\[
 D(e^{iv(xa+yb)})
 =\int_0^1\frac{1-u}{u}\bigl[\tr(e^{iuvX})-1\bigr]\dd u
\]
$n$ times at $v=0$. On the left, compact support permits
differentiation in the smooth test-function topology. On the right,
the differentiated integrand is bounded for real $v$ by
$u^{n-1}(1-u)\|X\|^n$. Thus
\[
 i^nD((xa+yb)^n)
 =i^n\tr(X^n)\int_0^1u^{n-1}(1-u)\dd u
 =\frac{i^n\tr(X^n)}{n(n+1)},
\]
which proves the lemma.
\end{proof}

For nonnegative integers $r,s$, let
\[
 W_{r,s}(A,B)=\binom{r+s}{r}^{-1}S_{r+s,s}(A,B).
\]
Thus $W_{r,s}$ is the average of all words with $r$ copies of $A$
and $s$ copies of $B$. Comparing coefficients in
\eqref{eq:directional-moments}, we obtain
\begin{equation}\label{eq:mixed-moments}
 D(a^rb^s)=\frac{\tr(W_{r,s}(A,B))}{(r+s)(r+s+1)}
 \qquad(r+s\geq1).
\end{equation}

\section{A scalar measure associated with a selfadjoint pair}\label{sec:scalar}

The following trace identity gives the sign needed to apply the
Herglotz representation theorem.

\begin{lemma}\label{lem:lyapunov}
Let $0\ne P\in\M$ be a projection and let $T\in P\M P$ satisfy
$\spec(T)\subseteq\{z\in\C:\im z<0\}$. Then
\[
 V=\int_0^\infty e^{itT^*}e^{-itT}\dd t
\]
converges in norm to a positive invertible element of $P\M P$, and
\begin{equation}\label{eq:lyapunov}
 -\im\tr(T)=\tfrac12\tr(V^{-1}).
\end{equation}
\end{lemma}
\begin{proof}
Choose a contour for the holomorphic functional calculus which encloses
$\spec(T)$ and lies in $\im z\leq-\delta$ for some $\delta>0$.
The contour formula for the exponential gives
$\|e^{-itT}\|\leq Ce^{-\delta t}$ for $t\geq0$. The integral therefore
converges in norm. Since its integrand is positive and converges to $P$ in
norm as $t\downarrow0$, we obtain $V\geq cP$ for some $c>0$.
Differentiating the integrand and integrating from zero to infinity,
we obtain
\[
 iT^*V-iVT=-P,\qquad \im(VT)=-\tfrac12P.
\]
The similar operator $Z=V^{1/2}TV^{-1/2}$ satisfies
\[
 \im Z=V^{-1/2}\im(VT)V^{-1/2}=-\tfrac12V^{-1}.
\]
Since $\tr(Z)=\tr(T)$, taking traces gives
\eqref{eq:lyapunov}.
\end{proof}

\begin{proposition}\label{prop:scalar}
Let $S,J\in\M$ be selfadjoint and put
$I_J=[\min\spec(J),\max\spec(J)]$. There is a unique finite positive
Borel measure $\nu_{S,J}$ supported in $I_J$ such that
\begin{equation}\label{eq:scalar-cauchy}
 \int\frac{\dd\nu_{S,J}(b)}{x-b}
 =\tr\left(\left[(x\one-J)^{-1/2}S(x\one-J)^{-1/2}\right]_+\right)
 \quad(x>\max\spec(J)).
\end{equation}
Its mass is $\tr(S_+)$. If $(S_n,J_n)\to(S,J)$ in operator norm,
with every pair selfadjoint, then $\nu_{S_n,J_n}\to\nu_{S,J}$ weakly.
\end{proposition}
\begin{proof}
Assume first that $S$ is invertible. For $z\in\C$ with $\im z>0$, put
\[
 X(z)=S(z\one-J)^{-1}.
\]
We claim that $\spec(X(z))\cap\R=\varnothing$. Indeed, $X(z)$ is
invertible, and for real $r\ne0$ we have
\[
 X(z)-r\one=[S-r(z\one-J)](z\one-J)^{-1}.
\]
The first factor has imaginary part $-r\im z\,\one$, which is
strictly definite. Indeed, for an operator $Y$ with $\im Y\geq c\one$ or
$\im Y\leq-c\one$, where $c>0$, the inequality
$\|Yv\|\|v\|\geq|\im\langle Yv,v\rangle|\geq c\|v\|^2$ gives
closed range and injectivity, and the same inequality for $Y^*$ gives
dense range. Thus the first factor is invertible, which proves the claim.

Let $E_-(z)$ be the Riesz idempotent selecting the part of
$\spec(X(z))$ in the lower half-plane. Near any $z_0$, a common contour
$\Gamma_-$ suffices: choose neighborhoods of the two compact spectral
parts with closures in their respective half-planes. Upper semicontinuity
keeps nearby spectra in their union. Then
\[
 E_-(z)=\frac1{2\pi i}\int_{\Gamma_-}(w\one-X(z))^{-1}\dd w
\]
\begingroup\postdisplaypenalty=10000
near $z_0$. Thus $E_-$ is locally holomorphic. Since the local expressions
agree on overlaps,
\begin{equation}\label{eq:scalar-holomorphic}
 m(z)=\tr(X(z)E_-(z))
\end{equation}
is holomorphic in the upper half-plane.
\par\endgroup

To determine the imaginary part, let $P$ be the orthogonal projection
onto the range of $E_-(z)$.
This range is closed and $P$, the support projection of
$E_-(z)E_-(z)^*$, belongs to $\M$. Since $X(z)$ commutes with
$E_-(z)$, the range of $E_-(z)$ is invariant under $X(z)$.
If $P\ne0$, the Riesz decomposition shows that the restriction
$T=PX(z)P$ to $P\mathcal H$
has spectrum equal to the selected lower spectral part. The faithful
representation of $P\M P$ on $P\mathcal H$ is inverse closed, so this
is also the spectrum in the corner. Relative to $P\oplus(\one-P)$,
we can write
\[
 E_-(z)=\begin{pmatrix}P&L\\0&0\end{pmatrix},\qquad
 X(z)=\begin{pmatrix}T&C\\0&R\end{pmatrix},\qquad
 X(z)E_-(z)=\begin{pmatrix}T&TL\\0&0\end{pmatrix}.
\]
Since the trace of an off-diagonal corner is zero,
$m(z)=\tr(T)$. Lemma~\ref{lem:lyapunov} therefore gives $\im m(z)\leq0$.
For $P=0$, we have $m(z)=0$, so the same inequality holds.

As $y\to\infty$, the resolvent expansion gives
$yX(iy)\to-iS$ in norm. The lower Riesz projection of $-iS$ is
$\one_{(0,\infty)}(S)$. Continuity of the corresponding contour
integrals therefore gives
\begin{equation}\label{eq:mass-limit}
 y m(iy)\longrightarrow-i\tr(S_+).
\end{equation}
Here multiplication by the positive scalar $y$ preserves the lower
spectral part.

To locate the support of the representing measure, we identify the real
boundary values. Write
$j_- =\min\spec(J)$ and $j_+=\max\spec(J)$, and fix $x>j_+$.
The operator $X(x)$ is similar to the invertible selfadjoint operator
\[
 K_x=(x\one-J)^{-1/2}S(x\one-J)^{-1/2}.
\]
Since its positive and negative spectra are separated, a contour around
the positive spectrum defines a holomorphic spectral idempotent near $x$.
To identify it with $E_-(z)$ above the real axis, let
$z=x'+iy$, where $x'>j_+$ and $y>0$, and let
$w=u+iv\in\spec(X(z))$. Since
\[
 \im[S-w(z\one-J)]=-[v(x'\one-J)+uy\one],
\]
the first factor in
$X(z)-w\one=[S-w(z\one-J)](z\one-J)^{-1}$ would be invertible if
$u>0$, $v\geq0$, or if $u<0$, $v\leq0$.
The case $u=0$, $v\ne0$ is excluded in the same way.
Thus the positive spectral part at $x$ continues into the lower
half-plane, and the negative part into the upper half-plane.
The two idempotents therefore agree, and similarity and traciality give
\begin{equation}\label{eq:right-boundary}
 m(x)=\tr(K_x\one_{(0,\infty)}(K_x))=\tr((K_x)_+).
\end{equation}
In particular, $m$ extends holomorphically across $x$ and is real there.
For $x<j_-$, the same calculation with $x'\one-J<0$ selects the
negative real spectral part of $X(x)$. Hence $m$ again extends
holomorphically and takes the real value
\[
 m(x)=-\tr\left(\left[(J-x\one)^{-1/2}S(J-x\one)^{-1/2}\right]_+\right).
\]

The scalar Herglotz theorem applied to $-m$ (see
\cite[Theorem~3.20]{Teschl}) gives
\[
 m(z)=c-dz+\int_{\R}
          \left(\frac1{z-b}+\frac b{1+b^2}\right)\dd\nu(b),
\]
where $c\in\R$, $d\geq0$, and $\nu$ is positive with
$\int(1+b^2)^{-1}\dd\nu(b)<\infty$. Taking imaginary parts at $iy$
gives
\[
 -y\im m(iy)=dy^2+\int\frac{y^2}{y^2+b^2}\dd\nu(b).
\]
It follows from \eqref{eq:mass-limit} and monotone convergence that
$d=0$ and $\nu(\R)=\tr(S_+)$. We can therefore absorb the correction
term in a real constant and write
$m(z)=c'+\int(z-b)^{-1}\dd\nu(b)$.
Since $m(iy)\to0$ and the integral is bounded by $\nu(\R)/y$, we
obtain $c'=0$.

To determine the support, take an interval $[u,v]$ contained in
$\R\setminus I_J$. Real holomorphic continuation implies that
$\im m(x+i\varepsilon)\to0$ uniformly on $[u,v]$. On the other hand,
\[
 -\frac1\pi\int_u^v\im m(x+i\varepsilon)\dd x
 =\int\frac1\pi\left[\arctan\frac{v-b}{\varepsilon}
                  -\arctan\frac{u-b}{\varepsilon}\right]\dd\nu(b).
\]
The integrand on the right lies in $[0,1]$ and tends to one for
$u<b<v$. Since the left side tends to zero and the endpoint contributions
are nonnegative, dominated convergence gives $\nu((u,v))=0$.
Such intervals cover $\R\setminus I_J$, and hence $\supp\nu\subseteq I_J$.
Together with \eqref{eq:right-boundary}, this proves existence for
invertible $S$. If two finite measures supported in $I_J$ satisfy
\eqref{eq:scalar-cauchy}, expansion at large positive $x$ shows that
all their moments agree. Polynomial density in $C(I_J)$ then proves
uniqueness.

For general selfadjoint $S$, choose invertible selfadjoint $S_n$ with
$\|S_n-S\|\to0$ by bounded Borel functional calculus. For example,
replace spectral values in $[0,1/n)$ by $1/n$ and values in $(-1/n,0)$
by $-1/n$. The measures $\nu_{S_n,J}$ have common compact support
$I_J$ and masses $\tr((S_n)_+)\to\tr(S_+)$. They therefore have weakly
convergent subsequences. The right side of
\eqref{eq:scalar-cauchy} converges by norm continuity of continuous
functional calculus on bounded selfadjoint sets. Thus every weak limit
satisfies \eqref{eq:scalar-cauchy}, so uniqueness identifies all these
limits. This defines $\nu_{S,J}$ and proves the mass formula.

Finally, suppose $(S_n,J_n)\to(S,J)$ in norm. The intervals $I_{J_n}$
are contained in one compact interval, their endpoints converge to
those of $I_J$, and the masses are bounded. Every weak cluster point
is therefore supported in $I_J$. For fixed $x$ larger than all the
right endpoints, the right side of \eqref{eq:scalar-cauchy} converges
to its value for $S,J$. Uniqueness identifies every cluster point
with $\nu_{S,J}$, and weak compactness gives convergence of the
whole sequence.
\end{proof}

For the moment calculation, substitute $x=t^{-1}$ in
\eqref{eq:scalar-cauchy}. For sufficiently small $t>0$, we obtain
\begin{equation}\label{eq:scalar-generating}
 \int\frac{\dd\nu_{S,J}(b)}{1-tb}
 =\tr\left(\left[(\one-tJ)^{-1/2}S(\one-tJ)^{-1/2}\right]_+\right).
\end{equation}
At $t=0$, the same equality follows from the mass formula.

\section{Weighted moment inequalities}\label{sec:positivity}

To compute the moments of the scalar measures, we use invariance of the
trace of a spectral projection under congruence.

\begin{lemma}\label{lem:inertia}
Let $H=H^*\in\M$ and let $C\in\M$ be invertible. Then
\[
 \tr(\one_{(0,\infty)}(C^*HC))
 =\tr(\one_{(0,\infty)}(H)).
\]
\end{lemma}
\begin{proof}
For projections $P,Q\in\M$, the polar decomposition of
$(\one-Q)P$ has initial projection $P-(P\wedge Q)$ and final
projection bounded by $\one-Q$. Taking traces gives
\[
 \tr(P\wedge Q)\geq\tr(P)+\tr(Q)-1.
\]
For $\varepsilon>0$, set $E_\varepsilon=\one_{[\varepsilon,\infty)}(H)$,
and let $Q_\varepsilon$ be the orthogonal projection onto the closed
subspace $C^{-1}E_\varepsilon\mathcal H$.
The polar decomposition of $C^{-1}E_\varepsilon$ shows that
$Q_\varepsilon\in\M$ and
$\tr(Q_\varepsilon)=\tr(E_\varepsilon)$. For a nonzero vector $v$
in this subspace, we have
\[
 \langle C^*HCv,v\rangle
 =\langle HCv,Cv\rangle\geq\varepsilon\|Cv\|^2>0.
\]
Its intersection with the nonpositive spectral subspace of $C^*HC$
is therefore zero. Applying the projection inequality, we obtain
$\tr(E_\varepsilon)\leq\tr(\one_{(0,\infty)}(C^*HC))$.
Letting $\varepsilon\downarrow0$ and using normality of the trace
proves one inequality. Apply it to $C^*HC$ and $C^{-1}$ to obtain
the reverse inequality.
\end{proof}

For selfadjoint $H\in\M$, the spectral theorem and Tonelli's theorem
give the layer-cake formula
\begin{equation}\label{eq:layer-cake}
 \tr(H_+)=\int_0^\infty\tr(\one_{(q,\infty)}(H))\dd q.
\end{equation}
For $H\geq0$ and an integer $r\geq1$, we similarly have
\begin{equation}\label{eq:power-layer-cake}
 \tr(H^r)=r\int_0^\infty q^{r-1}
                         \tr(\one_{(q,\infty)}(H))\dd q.
\end{equation}
Indeed, these follow by integrating the scalar identities for $x_+$
and $x^r$ against the spectral distribution of $H$.

\begin{proposition}\label{prop:weighted}
For positive $A,B\in\M$, let $D$ be the distribution from
Section~\ref{sec:fourier}. Then
\[
 D(a|p|^2)\geq0,\qquad D(b|p|^2)\geq0
 \quad\text{for every }p\in\C[a,b].
\]
\end{proposition}
\begin{proof}
Assume first that $A$ is invertible. Put
$\alpha=\|A\|$, $\beta=\|B\|$,
$\Theta=A^{-1/2}BA^{-1/2}$ and $L=\|\Theta\|$.
For $0<a<\alpha$, apply Proposition~\ref{prop:scalar} to
\[
 S_a=a^{-1}\one-A^{-1},\qquad J_a=a\Theta,
 \qquad \nu_a=\nu_{S_a,J_a}.
\]
Then $\nu_a$ is positive, supported in $[0,aL]$, and
\begin{equation}\label{eq:slice-bound}
 \nu_a(\R)=\tr((S_a)_+)\leq a^{-1}.
\end{equation}
The weak continuity in Proposition~\ref{prop:scalar} shows that
the moments
$h_s(a)=\int b^s\dd\nu_a(b)$ are continuous on $(0,\alpha)$.
Here the supports lie in $[0,\alpha L]$, and the masses are bounded
on each compact subinterval of $(0,\alpha)$.

To compute the generating function of these moments, choose $t\geq0$
with $t\alpha L<1$ and $t\beta<1$, and write
\[
 F_a(t)=\int\frac{\dd\nu_a(b)}{1-tb}.
\]
By \eqref{eq:scalar-generating}, we have
\begin{equation}\label{eq:slice-positive-part}
 F_a(t)=\tr((K_a(t))_+),\qquad
 K_a(t)=(\one-at\Theta)^{-1/2}S_a(\one-at\Theta)^{-1/2}.
\end{equation}
Set $R_a=\one-at\Theta$. Since
\[
 a^{-1}R_a-S_a
 =A^{-1}-t\Theta=A^{-1/2}(\one-tB)A^{-1/2}\geq0,
\]
we obtain $K_a(t)\leq a^{-1}\one$. It follows from
\eqref{eq:layer-cake} that
\[
 F_a(t)=\int_0^{1/a}
             \tr(\one_{(0,\infty)}(K_a(t)-q\one))\dd q.
\]
We apply congruence by $R_a^{1/2}$, which replaces
$K_a(t)-q\one$ by $S_a-qR_a$.
A second congruence, by $A^{1/2}$, followed by multiplication by $a$,
gives
\[
 aA^{1/2}(S_a-qR_a)A^{1/2}
 =(1-aq)A-a\one+a^2qtB.
\]
Set $v=aq$. Then the last expression is
$(1-v)A-a(\one-vtB)$, whose congruence by
$(\one-vtB)^{-1/2}$ is $M_v(t)-a\one$, where
\[
 M_v(t)=(1-v)(\one-vtB)^{-1/2}A(\one-vtB)^{-1/2}.
\]
Lemma~\ref{lem:inertia} applies at each step. Substituting
$\dd q=\dd v/a$, we obtain
\begin{equation}\label{eq:slice-layer-cake}
 F_a(t)=\frac1a\int_0^1\tr(\one_{(a,\infty)}(M_v(t)))\dd v.
\end{equation}
For $0\leq v\leq1$, we have
\[
 0\leq M_v(t)\leq\alpha(1-v)(\one-vtB)^{-1}
 \leq\alpha\frac{1-v}{1-vt\beta}\one\leq\alpha\one.
\]

Let $r\geq1$ be an integer. The bound
$a^rF_a(t)\leq a^{r-1}/(1-t\alpha L)$ proves integrability in $a$.
The spectral-tail integrand in \eqref{eq:slice-layer-cake} is jointly
Borel measurable in $(a,v)$. Indeed, the functions
$f_n(x)=\min\{1,n\max\{x,0\}\}$ increase to
$\one_{(0,\infty)}(x)$, and
$(a,v)\mapsto\tr(f_n(M_v(t)-a\one))$ is continuous.
Normality of the trace gives the desired integrand as their
pointwise limit. We can therefore apply Tonelli's theorem and
\eqref{eq:power-layer-cake} to obtain
\begin{align}
 \int_0^\alpha a^r F_a(t)\dd a
 &=\int_0^1\int_0^\alpha a^{r-1}
             \tr(\one_{(a,\infty)}(M_v(t)))\dd a\dd v\notag\\
 &=\frac1r\int_0^1\tr(M_v(t)^r)\dd v\notag\\
 &=\frac1r\int_0^1(1-v)^r
               \tr\bigl((A(\one-vtB)^{-1})^r\bigr)\dd v.
 \label{eq:integrated-generating}
\end{align}
The last equality uses cyclicity of the trace.

To compare the coefficients of $t^s$ in
\eqref{eq:integrated-generating}, we first count the relevant words.
For $r\geq1$, $s\geq0$, put
\[
 R_{r,s}=\sum_{i_1+\cdots+i_r=s}
             \tr(AB^{i_1}\cdots AB^{i_r}),
 \qquad i_1,\ldots,i_r\geq0.
\]
Mark one occurrence of $A$ in a word with $r$ copies of $A$ and $s$
copies of $B$. The sum of the traces over these marked words is
$r\tr(S_{r+s,s}(A,B))$. Rotate each word to start at its marked
letter. It then has the form $AB^{i_1}\cdots AB^{i_r}$.
The weak composition $(i_1,\ldots,i_r)$, together with the original
position of the marked letter among the $r+s$ positions, recovers
the marked word. Cyclicity of the trace thus gives
\begin{equation}\label{eq:word-count}
 R_{r,s}=\frac r{r+s}\tr(S_{r+s,s}(A,B)).
\end{equation}
This count includes all marked positions, also in periodic words.

For $|t|\leq t_0$ with $t_0\alpha L<1$ and $t_0\beta<1$, the
geometric series on the left of \eqref{eq:integrated-generating}
is bounded after multiplication by $a^r$ by
$a^{r-1}/(1-t_0\alpha L)$. On the right, the resolvent series and
its finite products converge uniformly in $v\in[0,1]$.
These bounds justify expansion of both sides in convergent power series
near zero and termwise integration of their coefficients.
Since the series agree for small positive $t$, their coefficients agree.
For the coefficient of $t^s$, we obtain
\begin{align}
 \int_0^\alpha a^r h_s(a)\dd a
 &=\frac{R_{r,s}}r\int_0^1v^s(1-v)^r\dd v\notag\\
 &=\frac{R_{r,s}}r\frac{r!s!}{(r+s+1)!}
 =\frac{\tr(W_{r,s}(A,B))}{(r+s)(r+s+1)}.
 \label{eq:slice-moments}
\end{align}
In the second equality, the beta integral follows by repeated
integration by parts; the last equality uses \eqref{eq:word-count}.
When $L=0$, the same calculation applies with all positive-degree
$b$-moments equal to zero.

To deduce positivity, let $p(a,b)=\sum_{j,k}c_{jk}a^jb^k$ be a polynomial.
Each monomial
in $a|p|^2$ has positive $a$-degree. By \eqref{eq:mixed-moments}
and \eqref{eq:slice-moments}, we obtain
\begin{equation}\label{eq:gram}
 D(a|p|^2)
 =\int_0^\alpha a\int |p(a,b)|^2\dd\nu_a(b)\dd a\geq0.
\end{equation}
The integral is finite, since $p$ is bounded on
$[0,\alpha]\times[0,\alpha L]$ and $a\nu_a(\R)\leq1$.

If both $A$ and $B$ are invertible, interchanging the operators and
the two scalar coordinates gives $D(b|p|^2)\geq0$ as well.
We pass to general positive $A,B$ by applying these inequalities to
$A_\varepsilon=A+\varepsilon\one$ and
$B_\varepsilon=B+\varepsilon\one$.
For a fixed polynomial $p$, each of the two distribution pairings
is, by \eqref{eq:mixed-moments}, a finite linear combination of
traces of words in $A_\varepsilon,B_\varepsilon$. The words converge
in operator norm as $\varepsilon\downarrow0$, and hence their traces
converge to those for $A,B$. The limiting inequalities prove the
proposition.
\end{proof}

\section{The planar measure}\label{sec:measure}

To obtain positive measures from the weighted inequalities, we pass
from polynomial squares to nonnegative smooth functions.

\begin{lemma}\label{lem:positive-distribution}
Let $T$ be a compactly supported distribution on $\R^2$ such that
$T(|p|^2)\geq0$ for every $p\in\C[a,b]$. Then $T$ is a finite
positive Borel measure.
\end{lemma}
\begin{proof}
Choose a closed cube $Q$ whose interior contains $\supp T$, and a
smooth compactly supported cutoff in its interior which is one near
$\supp T$. Since a distribution has finite order on a fixed compact
set, there are an integer $m$ and a constant $C$ such that
\[
 |T(f)|\leq C\max_{|\gamma|\leq m}\sup_Q|\partial^\gamma f|
 \qquad(f\in C^\infty(\R^2)).
\]
To approximate in this $C^m$ norm, let $h$ be real and smooth near
$Q$. We multiply $h$ by a smooth cutoff
equal to one near $Q$ and supported in its domain. Convolving the
resulting compactly supported function with centered Gaussians gives
functions converging to $h$ in $C^m(Q)$. Each convolution extends to
an entire function on $\C^2$, since its defining integral has compact
support. Its Taylor polynomials have real coefficients and
converge with all derivatives uniformly on $Q$. A diagonal choice
therefore gives real polynomials approximating $h$ in $C^m(Q)$.

Let $\varphi\in C_c^\infty(\R^2)$ be nonnegative and let
$\varepsilon>0$. We apply this approximation to
$h=(\varphi+\varepsilon)^{1/2}$. The product rule shows that the
squares of the approximating polynomials converge to
$\varphi+\varepsilon$ in $C^m(Q)$. Hence
\[
 T(\varphi+\varepsilon)=\lim_n T(p_n^2)\geq0.
\]
Letting $\varepsilon\downarrow0$ gives $T(\varphi)\geq0$.
For a real smooth function $f$ supported in a compact set $L$, choose
$\psi\in C_c^\infty(\R^2)$ with $\psi\geq0$ and $\psi=1$ on $L$.
Then
$-\|f\|_\infty\psi\leq f\leq\|f\|_\infty\psi$, so
$|T(f)|\leq\|f\|_\infty T(\psi)$.
The functional therefore extends to continuous compactly supported
functions, and the Riesz representation theorem gives a positive Radon
measure.
Its support is compact, so its mass is finite.
\end{proof}

For the positive pair $A,B$ and its distribution $D$,
Proposition~\ref{prop:weighted} and
Lemma~\ref{lem:positive-distribution} show that
\[
 \eta_A=aD,\qquad \eta_B=bD
\]
are finite positive measures supported in $K$, with masses
\begin{equation}\label{eq:weighted-masses}
 \eta_A(K)=D(a)=\frac{\tr(A)}2,
 \qquad \eta_B(K)=D(b)=\frac{\tr(B)}2.
\end{equation}
To control their mass at the origin, we determine two coordinate
marginals.

Let $\sigma_A$ be the spectral distribution of $A$, characterized by
$\tr(f(A))=\int f(\lambda)\dd\sigma_A(\lambda)$, and for $x\geq0$ put
\[
 f_A(x)=\int_{\lambda>0}(1-x/\lambda)_+\dd\sigma_A(\lambda).
\]
This function takes values in $[0,1]$ and vanishes for $x>\alpha$.
For every integer $n\geq0$, Tonelli's theorem gives
\begin{align*}
 \int_0^\alpha x^n f_A(x)\dd x
 &=\int_{\lambda>0}\int_0^\lambda
                 x^n(1-x/\lambda)\dd x\dd\sigma_A(\lambda)\\
 &=\frac{\tr(A^{n+1})}{(n+1)(n+2)}
 =D(a^{n+1})=\int_K a^n\dd\eta_A(a,b).
\end{align*}
Since these two finite measures on $[0,\alpha]$ have the same moments,
polynomial density shows that they agree. Thus the first-coordinate
marginal of $\eta_A$ is $f_A(x)\dd x$. Interchanging $A$ and $B$
gives the second-coordinate marginal of $\eta_B$. It follows that
\begin{equation}\label{eq:no-axis-atoms}
 \eta_A(\{a=0\})=0,\qquad \eta_B(\{b=0\})=0.
\end{equation}
The same conclusion holds when an operator is zero, since the
corresponding weighted measure then has mass zero.

We set $\omega=\eta_A+\eta_B$ and define a measure on $K\setminus\{0\}$ by
\begin{equation}\label{eq:divide}
 \dd\mu(a,b)=\frac1{a+b}\dd\omega(a,b).
\end{equation}
Since $a+b$ is positive on this set, $\mu$ is a positive Radon measure
on the punctured plane. By \eqref{eq:no-axis-atoms}, $\omega$ has no
atom at the origin. Hence \eqref{eq:weighted-masses} gives
\begin{equation}\label{eq:mu-first}
 \int(a+b)\dd\mu=\omega(K)=\frac{\tr(A+B)}2.
\end{equation}
Multiplication of distributions by the scalar coordinates commutes, so
$b\eta_A=a\eta_B$. Therefore
\[
 (a+b)\eta_A=a\omega,\qquad (a+b)\eta_B=b\omega.
\]
Dividing these identities away from the origin and using
\eqref{eq:no-axis-atoms}, we obtain
\begin{equation}\label{eq:separate-weights}
 a\mu=\eta_A,\qquad b\mu=\eta_B,
\end{equation}
where the weighted measures are extended by zero at the origin.

For every smooth $\varphi$, the fundamental theorem of calculus
along the segment from zero to $(a,b)\in K$ gives
\[
 |\varphi(a,b)-\varphi(0,0)|
 \leq(a+b)\max_{x\in K}\max\{|\partial_a\varphi(x)|,
                                  |\partial_b\varphi(x)|\}.
\]
Together with \eqref{eq:mu-first}, this shows that
\[
 D_\mu(\varphi)=\int[\varphi(a,b)-\varphi(0,0)]\dd\mu(a,b)
\]
defines a distribution supported in $K$.
By \eqref{eq:separate-weights}, we have
$aD_\mu=aD$, $bD_\mu=bD$, and $D_\mu(1)=D(1)=0$.
Let $R=D-D_\mu$. For any smooth $\varphi$, write
\[
 \varphi(a,b)-\varphi(0,0)
 =a\int_0^1\partial_a\varphi(ta,tb)\dd t
  +b\int_0^1\partial_b\varphi(ta,tb)\dd t.
\]
Applying the compactly supported distribution $R$ and using
$aR=bR=0$, $R(1)=0$, we obtain $R(\varphi)=0$. Thus
\begin{equation}\label{eq:compensated}
 D(\varphi)=\int[\varphi(a,b)-\varphi(0,0)]\dd\mu(a,b)
 \qquad(\varphi\in C^\infty(\R^2)).
\end{equation}
In particular, \eqref{eq:directional-moments} gives
\begin{equation}\label{eq:measure-moments}
 \int(xa+yb)^n\dd\mu(a,b)
 =\frac{\tr((xA+yB)^n)}{n(n+1)}
 \qquad(n\geq1).
\end{equation}
Every polynomial vanishing at zero is bounded on $K$ by a constant
times $a+b$, so these integrals are absolutely convergent.

\begin{proof}[Proof of Theorem~\ref{thm:planar}]
The positive measure $\mu$ constructed above satisfies
\eqref{eq:first-moment} and \eqref{eq:measure-moments}.
To obtain the identity for Borel functions, fix real $x,y$ and put
$\ell(a,b)=xa+yb$ and $X=xA+yB$.
Let $\sigma_X$ be the spectral distribution of $X$.
We define finite positive measures $\lambda$ and $\lambda_0$ on $\R$ by
\begin{align*}
 \lambda(E)&=\int\ell(a,b)^2\one_E(\ell(a,b))\dd\mu(a,b),\\
 \lambda_0(E)&=\int_\R\int_0^1
               (us)^2\one_E(us)\frac{1-u}{u}\dd u\dd\sigma_X(s).
\end{align*}
These measures have support in a common compact interval, and for
every $n\geq0$, \eqref{eq:measure-moments} gives
\[
 \int v^n\dd\lambda(v)
 =\frac{\tr(X^{n+2})}{(n+2)(n+3)}
 =\int v^n\dd\lambda_0(v).
\]
Polynomial density therefore gives $\lambda=\lambda_0$.
If $f:\R\to[0,\infty]$ is Borel with $f(0)=0$, integration of
$f(v)/v^2$ on $\R\setminus\{0\}$ yields
\begin{equation}\label{eq:borel-marginal}
 \int f(\ell(a,b))\dd\mu(a,b)
 =\int_\R\int_0^1\frac{1-u}{u}f(us)\dd u\dd\sigma_X(s).
\end{equation}
Both sides are understood as nonnegative integrals.
For a complex-valued $f$ satisfying the hypotheses of the theorem,
choose $R>\|X\|$. For $|s|\leq\|X\|$, substitution gives
\[
 \int_0^1\frac{1-u}{u}|f(us)|\dd u
 \leq\int_{-R}^{R}\frac{|f(v)|}{|v|}\dd v.
\]
At $s=0$, the left side is zero. We apply
\eqref{eq:borel-marginal} to $|f|$ to obtain absolute integrability,
and then to the positive and negative parts of its real and imaginary
parts. The spectral theorem gives \eqref{eq:planar-marginal}.

To prove uniqueness, observe that any measure satisfying the theorem
has the moments
\eqref{eq:measure-moments}. Expanding in $x,y$ determines all mixed
moments of positive degree. Thus the finite measure $(a+b)\mu$ has
prescribed integrals against every polynomial on $K$.
Polynomial density determines this measure uniquely, and division by
$a+b$ on $K\setminus\{0\}$ proves uniqueness of $\mu$.
\end{proof}

\section{Proof of Theorem A}\label{sec:bmv}

Let $A,B\in\M$ be positive and let $\mu=\mu_{A,B}$ be the measure
from Theorem~\ref{thm:planar}. For $m\geq1$, comparison of coefficients
in \eqref{eq:measure-moments} gives
\begin{equation}\label{eq:coefficient-positivity}
 \tr(S_{m,k}(A,B))
 =m(m+1)\binom mk\int a^{m-k}b^k\dd\mu(a,b)\geq0
 \qquad(0\leq k\leq m).
\end{equation}
The integrals are finite by \eqref{eq:mu-first} and boundedness of
$K$. For $m=0$, the assertion is $\tr(\one)=1$.

Let $H=H^*\in\M$, let $V\geq0$, and put
$f(z)=\tr(e^{H-zV})$ for $z\in\C$. The exponential series converges
locally uniformly in operator norm, so $f$ is entire.
Fix $t_0\geq0$ and choose $c\in\R$ such that
$C=H-t_0V+c\one\geq0$. Then
\begin{equation}\label{eq:shifted-exponential}
 f(t_0-s)=e^{-c}\tr(e^{C+sV})
 =e^{-c}\sum_{m=0}^\infty\frac1{m!}
                      \sum_{k=0}^m s^k\tr(S_{m,k}(C,V)).
\end{equation}
The estimate
\begin{equation}\label{eq:word-bound}
 \|S_{m,k}(C,V)\|\leq\binom mk\|C\|^{m-k}\|V\|^k
\end{equation}
bounds the sum of the norms by $e^{\|C\|+|s|\|V\|}$.
We can therefore collect coefficients and differentiate
\eqref{eq:shifted-exponential} term by term. By
\eqref{eq:coefficient-positivity}, we obtain
\begin{equation}\label{eq:derivative-signs}
 (-1)^kf^{(k)}(t_0)
 =e^{-c}k!\sum_{m=k}^\infty\frac{\tr(S_{m,k}(C,V))}{m!}\geq0.
\end{equation}
This proves complete monotonicity.

We construct the Laplace measure explicitly. For each integer $N\geq1$,
put
\[
 w_{N,k}=\frac{N^k(-1)^kf^{(k)}(N)}{k!},\qquad
 \rho_N=\sum_{k=0}^\infty w_{N,k}\delta_{k/N}.
\]
The weights are nonnegative by \eqref{eq:derivative-signs}.
Taylor expansion of the entire function $f$ about $N$ gives
\begin{equation}\label{eq:discrete-generating}
 \sum_{k=0}^\infty w_{N,k}z^k=f(N-Nz)\qquad(z\in\C).
\end{equation}
At $z=1$ we obtain $\rho_N([0,\infty))=f(0)$. Differentiating
\eqref{eq:discrete-generating} and then putting $z=1$ gives
\[
 \int s\dd\rho_N(s)=\frac1N\sum_{k=0}^\infty k w_{N,k}=-f'(0).
\]
Consequently,
$\rho_N([R,\infty))\leq-f'(0)/R$ for $R>0$.
The measures therefore have fixed finite mass and are tight. Weak
compactness then gives a subsequence converging weakly to a finite
positive measure $\rho$ on $[0,\infty)$, of mass $f(0)$.
For $t\geq0$, substitution of $z=e^{-t/N}$ in
\eqref{eq:discrete-generating} yields
\[
 \int e^{-ts}\dd\rho_N(s)=f\bigl(N(1-e^{-t/N})\bigr).
\]
Since the integrand is bounded and continuous and
$N(1-e^{-t/N})\to t$, we can pass to the limit and obtain
\begin{equation}\label{eq:laplace}
 f(t)=\int e^{-ts}\dd\rho(s)\qquad(t\geq0).
\end{equation}

If two finite measures have this Laplace transform, push them forward
by $s\mapsto e^{-s}$ and regard them as measures on $[0,1]$ with
zero mass at zero. Their integer moments agree by
\eqref{eq:laplace}, including the zeroth moment. Polynomial density
on $[0,1]$ gives equality of the pushforwards and hence of the original
measures, which proves uniqueness.

It remains to determine the support. For $t>0$, differentiation
under the integral in \eqref{eq:laplace} is justified by boundedness
of $s^ne^{-ts}$ on $[0,\infty)$, locally uniformly in $t>0$.
Letting $t\downarrow0$, monotone convergence gives
\begin{equation}\label{eq:laplace-moments}
 \int s^n\dd\rho(s)=(-1)^nf^{(n)}(0)\qquad(n\geq0).
\end{equation}
The same word estimate as in \eqref{eq:word-bound}, applied to $H,V$,
gives
\begin{align*}
 |f^{(n)}(0)|
 &\leq n!\sum_{m=n}^\infty
       \frac{\binom mn\|H\|^{m-n}\|V\|^n}{m!}\\
 &=e^{\|H\|}\|V\|^n.
\end{align*}
If $\rho([R,\infty))=d>0$ for some $R>\|V\|$, then
\eqref{eq:laplace-moments} would imply
$dR^n\leq e^{\|H\|}\|V\|^n$ for every $n$, which is impossible.
Thus $\supp\rho\subseteq[0,\|V\|]$.

Write $v_- =\min\spec(V)$. Since
$H-tV\leq(\|H\|-tv_-)\one$ for $t\geq0$, spectral calculus gives
\[
 f(t)\leq e^{\|H\|-tv_-}.
\]
If $\rho([0,v_--\delta])=d>0$ for some $\delta>0$, then
\eqref{eq:laplace} would give
$d e^{-t(v_--\delta)}\leq e^{\|H\|-tv_-}$ for all $t\geq0$,
again a contradiction. It follows that
$\supp\rho\subseteq[v_-,\|V\|]$.
As $V\geq0$, we have $\|V\|=\max\spec(V)$, which proves the
support assertion. Finally, $\rho(\R)=f(0)=\tr(e^H)$, and
Theorem~\ref{thm:bmv} follows.

\section*{Acknowledgments}
OpenAI's language models were
used in developing and checking the arguments and in preparing the
manuscript.

\end{document}